\documentclass[preprint,12pt]{elsarticle}
\usepackage{times}
\usepackage{amsfonts,amsmath,amssymb,amsthm}
\usepackage{color,soul}
\biboptions{numbers,sort&compress} % options to natbib (included in elsarticle)
\let\ts=\textstyle
\def\fracskip{\mskip 1mu \relax}
\let\oldfrac=\frac
\def\nfrac#1#2{\oldfrac{\fracskip#1\fracskip}{\fracskip#2\fracskip}}
\def\tfrac#1#2{{\ts\nfrac{#1}{#2}}}
\let\frac=\nfrac
\def\pd#1#2{\frac{\partial#1}{\partial#2}}
\def\pdd#1#2#3{\ifx#2#3\pd{^2#1}{#2^2}\else\pd{^2#1}{#2\partial#3}\fi }

\let\BL=\biggl \let\BR=\biggr
\let\bl=\bigl \let\br=\bigr

\def\arbs{are arbitrary constants}

\def\Equation#1. {\medbreak{\bfseries\itshape{Equation\kern.3333em\relax#1.}}\enspace\ignorespaces }

\newtheoremstyle{remark}{\medskipamount}{\medskipamount}
  {\small\rmfamily}{\parindent}{\footnotesize\sffamily}{.}{.5em}{}
\theoremstyle{remark}
\newtheorem{remark}{Remark}%[section]

\let\ts=\textstyle
\let\bl=\bigl \let\br=\bigr
 
\let\BL=\biggl \let\BR=\biggr

\begin{document}
\begin{frontmatter}
\title{Hypersingular nonlinear boundary-value problems\\ with a small parameter$^{***}$}
\author[ipm,bmstu,mephi]{Andrei D. Polyanin\corref{cor1}}
  \ead{polyanin@ipmnet.ru}
\author[us]{Inna K. Shingareva\corref{cor2}}
  \ead{inna@mat.uson.mx}
  \cortext[cor1]{Principal corresponding author}
  \cortext[cor2]{Corresponding author}
\address[ipm]{Institute for Problems in Mechanics, Russian Academy of Sciences,\\
  101 Vernadsky Avenue, bldg~1, 119526 Moscow, Russia}
\address[bmstu]{Bauman Moscow State Technical University,\\
  5 Second Baumanskaya Street, 105005 Moscow, Russia}
\address[mephi]{National Research Nuclear University MEPhI,
  31 Kashirskoe Shosse, 115409 Moscow, Russia}
\address[us]{University of Sonora, Blvd. Luis Encinas y Rosales S/N, Hermosillo C.P. 83000, Sonora, M\'exico}

\begin{abstract}
For the first time, some hypersingular nonlinear boundary-value problems with a small parameter~$\varepsilon$ at the highest derivative are described. These
problems essentially (qualitatively and quantitatively) differ from the usual linear and quasilinear singularly perturbed boundary-value problems and have the
following unusual properties:\let\thefootnote\relax\footnotetext{
\hskip-20pt $^{***}$ This article will be published in the \textit{Applied Mathematics Letters\/} (2018).}

(i) in hypersingular boundary-value problems, super thin boundary layers arise,
and the derivative at the boundary layer can have very large values of the order
of $e^{1/\varepsilon}$ and more (in standard problems with boundary layers,
the derivative at the boundary has the order of $\varepsilon^{-1}$ or less);

(ii) in hypersingular boundary-value problems, the position of the boundary layer
is determined by the values of the unknown function at the boundaries (in standard problems
with boundary layers, the position of the boundary layer is determined
by coefficients of the given equation, and the values
of the unknown function at the boundaries do not play a role here);

(iii) hypersingular boundary-value problems do not admit a direct application
of the method of matched asymptotic expansions (without a preliminary nonlinear
point transformation of the equation under consideration).

Examples of hypersingular nonlinear boundary-value problems
with ODEs and PDEs are given and their exact solutions are obtained.
It is important to note that the exact solutions presented in this paper 
can be used to compare the effectiveness of various methods of numerical integration of
singularly perturbed problems with boundary layers, 
and also to develop new numerical and approximate analytical methods.
%It is shown that for some boundary values there can be
%no boundary layers as $\varepsilon\to 0$.
\end{abstract}

\end{frontmatter}

%%%%%%%%%%%%%%%%%%%% 1
\section{Introduction}\label{s:1}

Singularly perturbed boundary-value problems with a small parameter at the highest derivative
are often encountered in hydro- and aerodynamics, theory of mass and heat transfer,
theory of elasticity, nonlinear mechanics and other applications.
An important qualitative feature of singular boundary-value problems is that
for the zero value of a small parameter the order of the differential equation
under consideration decreases and some parts of the boundary conditions cannot be satisfied.
Various problems and solution methods for ODEs and PDEs
with a small parameter at the highest derivative are described, for example, in
\cite{eck1973,lag1988,mal1991,ker1996,nay2000,mir2001,Pol_Zai2003,ver2005,shish2009,kop2010,geng2016,Pol_Zai2017,Pol_Shing2018b}.

For singularly perturbed quasilinear boundary-value problems of the form
\begin{align}
\varepsilon y^{\prime\prime}_{xx}+py^{\prime}_x+q(x,y)=0\quad \ (0<x<1);\quad \ y(0)=\alpha,\quad \ y(1)=\beta
\label{01}
\end{align}
with a small parameter $\varepsilon$, the position of the boundary layer is determined
by the sign of the coefficient $p$.
For $p>0$, the boundary layer is formed at the left boundary near the point  $x=0$,
for $p<0$, at the right boundary near the point $x=1$.
An analogous situation holds for more complex singularly perturbed quasilinear boundary-value problems of the form \eqref{01} with $p=p(x)$.

We note that for $p>0$ and $\varepsilon\to 0$ the derivative of the solution of the problem~\eqref{01}
takes large values proportional to $\varepsilon^{-1}$ on the left boundary.

In this article we will consider hypersingular boundary-value problems with
a small parameter, the solutions of which qualitatively and quantitatively differ
from the solutions of the problems \eqref{01}.

%%%%%%%%%%%%%%%%%%%%%%%%%% 2
\section{Hypersingular boundary-value problems for ordinary differential\\ equations with a small parameter}\label{s:2}
%\section{A class of hypersingular nonlinear boundary-value problems with a small parameter}

\subsection{Example of hypersingular boundary-value problem. Exact solution,\\ qualitative features}\label{ss:2.1}

\textit{Problem 1.}
Consider the nonlinear boundary-value problem
\begin{align}
\varepsilon y^{\prime\prime}_{xx}+p(y^\prime_x)^2+qy^\prime_x=0;\label{e01}\\
y(0)=\alpha,\quad \ y(1)=\beta.\label{e01a}
\end{align}
In what follows, we assume that $p>0$, $q>0$ and $\varepsilon>0$ is a small parameter.
\goodbreak

Equation~\eqref{e01} admits an exact linearization by means of the transformation
$y=(\varepsilon/p)\ln u$ (see also Section~\ref{ss:2.2}). Its general solution is determined by the formula
\begin{align}
y=\frac{\varepsilon}p\ln(C_1+C_2e^{-qx/\varepsilon}), \label{e02}
\end{align}
where $C_1$ and $C_2$ \arbs.

The exact solution of the boundary-value problem~\eqref{e01}--\eqref{e01a} has the form
\begin{align}
y=\frac{\varepsilon}p\ln\BL[\frac{(e^{\alpha p/\varepsilon}
-e^{\beta p/\varepsilon})e^{-qx/\varepsilon}+e^{\beta p/\varepsilon}
-e^{(\alpha p-q)/\varepsilon}}{1-e^{-q/\varepsilon}}\BR].
\label{e03}
\end{align}
%\blh{Ïðîâåðåíî!}

Let us find the derivatives on the left and right boundaries:
\begin{align}
y^\prime_x|_{x=0}=\frac qp\frac{e^{(\beta-\alpha)p/\varepsilon}-1}{1-e^{-q/\varepsilon}},\quad \
y^\prime_x|_{x=1}=\frac qp\frac{1-e^{(\alpha-\beta)p/\varepsilon}}{e^{q/\varepsilon}-1}.
\label{e04}
\end{align}

The analysis shows that, depending on the values of the parameters $\alpha$ and $\beta$
(as \text{$\varepsilon\to 0$}), the following three qualitatively different situations are possible:
\begin{align}
\text{(i)}&\ \beta>\alpha\qquad\qquad \ \ \ \,\, (\text{there is a boundary layer on the left, near $x{=}0$});\notag\\
\text{(ii)}&\ \beta<\alpha-(q/p)\quad\ \ \, (\text{there is a boundary layer on the right, near $x{=}1$});\label{e05}\\
\text{(iii)}&\ \alpha{-}(q/p)\le \beta\le \alpha\ (\text{there is no boundary layer}).\notag
\end{align}

On the line $\beta=\alpha$, the problem~\eqref{e01}--\eqref{e01a} has the trivial solution $y=\alpha$, and
on the line $\beta=\alpha-(q/p)$ it has the linear solution $y=\alpha-(q/p)x$
(both these solutions are degenerate and do not depend on $\varepsilon$).
It is seen from~\eqref{e05} that the position of the boundary layer (or its absence)
in the problem under consideration is determined by the boundary values $\alpha$ and $\beta$,
which is qualitatively different from the situation typical
for standard problems with a boundary layer (see Section~1),
where the boundary conditions do not affect the position of the boundary layer.

\begin{figure}
\centering
\vskip-2pc
{\includegraphics[scale=0.34]{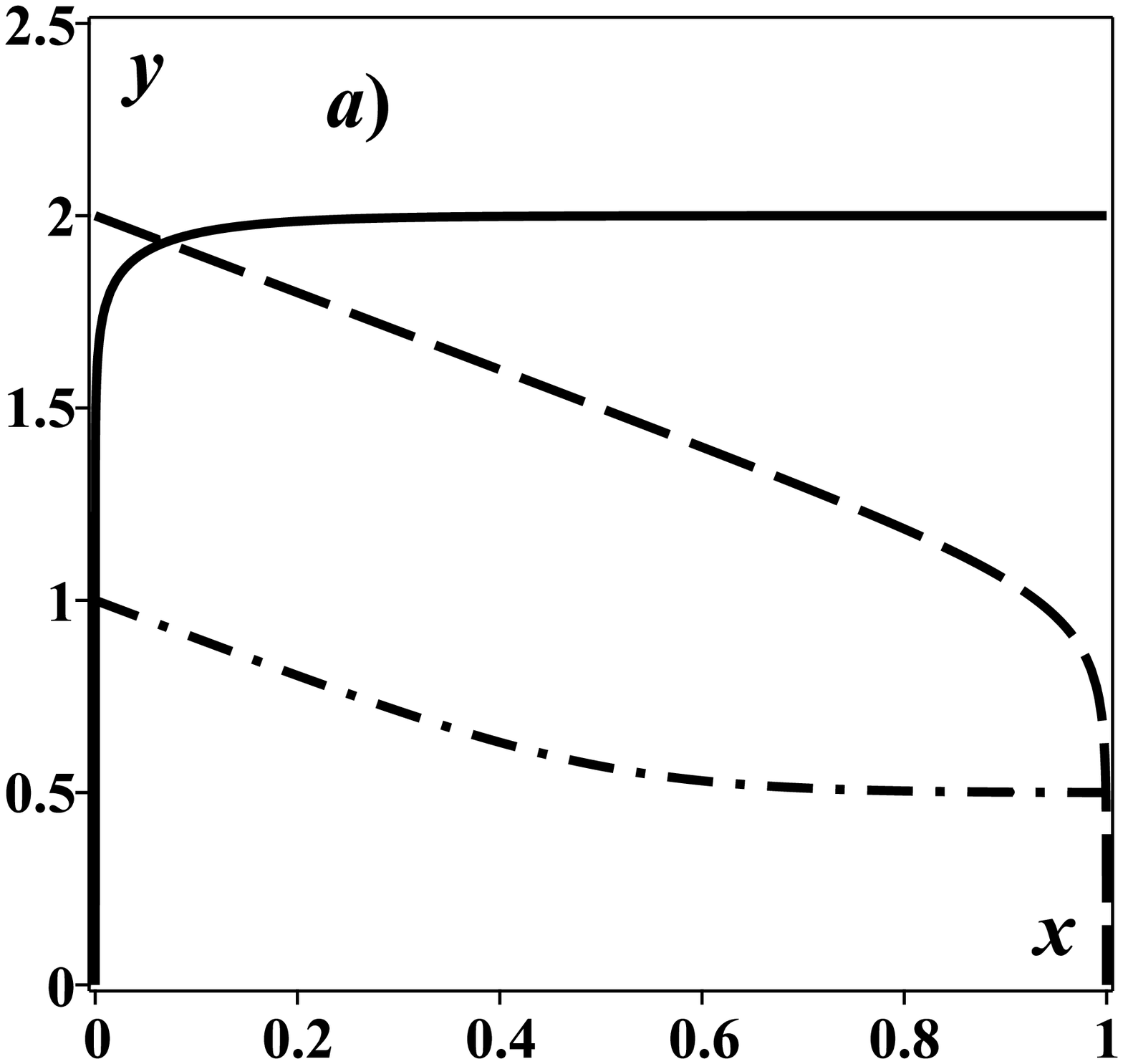}\ \includegraphics[scale=0.34]{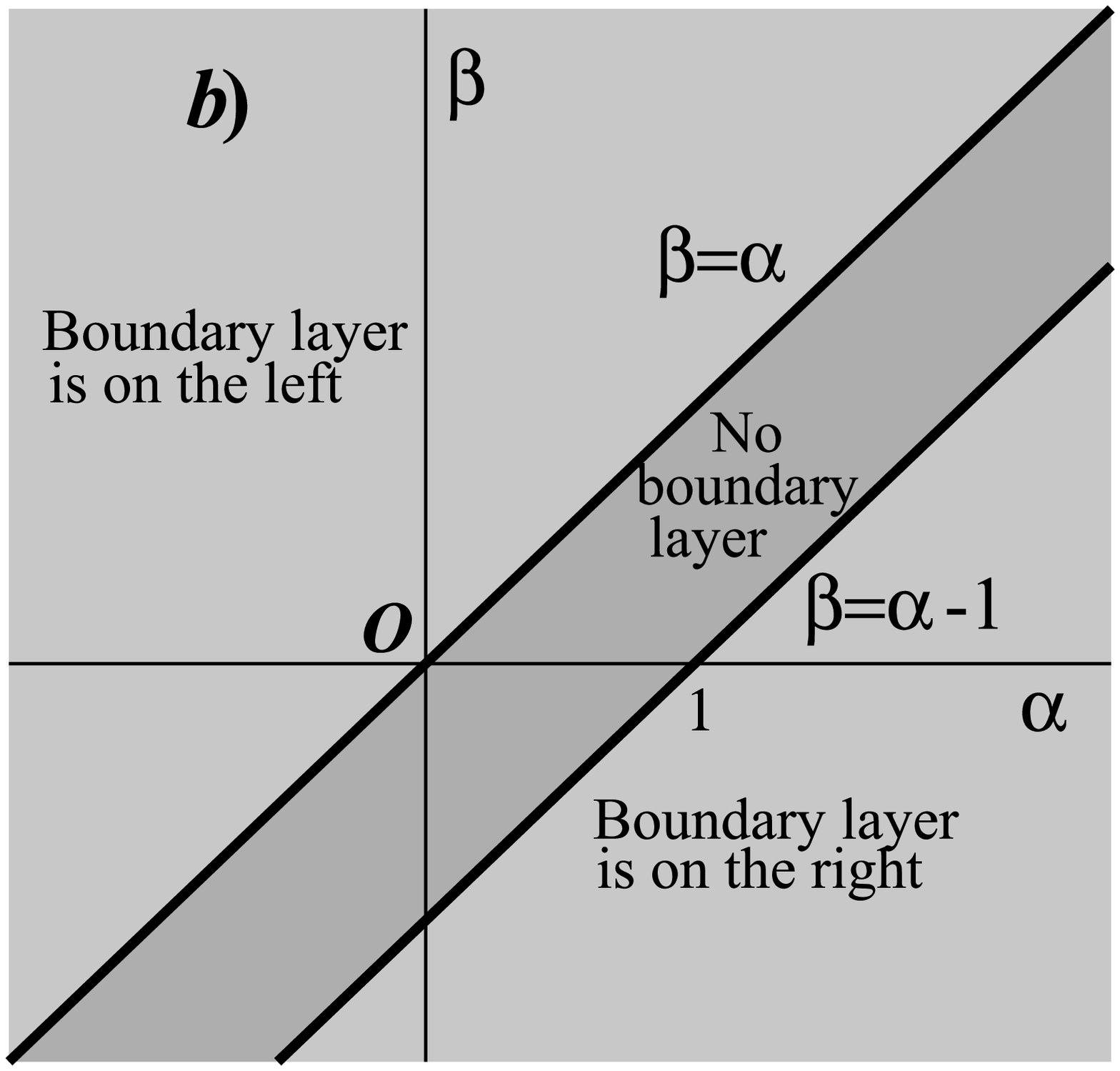}}
\vskip-1pc
\caption{a) Curves constructed using formula~\eqref{e03} for $p=q=1$ and $\varepsilon=0.1$
for  $\alpha=0$, $\beta=2$ (solid line); $\alpha=2$, $\beta=0$ (dashed line);
$\alpha=1$, $\beta=1/2$ (dash-dotted line); b) regions in the plane of the parameters
$\alpha$, $\beta$ with qualitatively different types of solutions of the problem~\eqref{e01}--\eqref{e01a}. }
\label{fig:Fig1}
\end{figure}

Fig.~\ref{fig:Fig1}\textit{a} represents the curves describing the exact solutions
of the problem~\eqref{e01}--\eqref{e01a} and corresponding to the three different cases in~\eqref{e05};
they are constructed according to the formula~\eqref{e03} for $p=q=1$ and $\varepsilon=0.1$
and are depicted by the solid line for $\alpha=1$, $\beta=1/2$; by the dashed line for $\alpha=2$, $\beta=0$;
by the dash-dotted line for $\alpha=1$, $\beta=1/2$.
Fig.~\ref{fig:Fig1}\textit{b} represents the plane of the parameters $\alpha$, $\beta$,
in which the parts of the plane defining three qualitatively different types of the solutions
of the problem~\eqref{e01}--\eqref{e01a} and satisfying the inequalities (i), (ii), (iii) in~\eqref{e05}
for $p=q=1$ are shaded in different ways; these parts of the plane are delimited by the
straight lines $\beta=\alpha$ and $\beta=\alpha-1$, which correspond to the
degenerate solutions.

%%%%%%%%%%%%%%%%%%%%  Figure 1

We set $p=q=1$, $\alpha=0$, $\beta=1$ in~\eqref{e01}--\eqref{e01a}. It follows from the first relation~\eqref{e04}
that even for moderately small $\varepsilon=\frac 1{10}$ the derivative reaches
extremely large values $y^\prime_x|_{x=0}\simeq e^{1/\varepsilon}\approx 22026.442$ on the left boundary,
which corresponds to the hyperfine (hypersingular) boundary layer.
For comparison, we note that in the problem~\eqref{e01}--\eqref{e01a} for $p=q=1$, $\alpha=0$, $\beta=1$, $\varepsilon=\frac 1{20}$
the derivative on the left boundary will be larger than the derivative
on the left boundary in the problem~\eqref{01} for $p=1$, $q(x,y)=y$, $\alpha=0$, $\beta=1$, $\varepsilon=10^{-8}$.

Similar problems differ fundamentally from the classical problems with a boundary layer
(in which the derivative at the boundary is of the order of $\varepsilon^{-1}$ or less)
and no such problems have been encountered in the literature (see, for example,
\cite{eck1973,lag1988,mal1991,ker1996,nay2000,mir2001,Pol_Zai2003,ver2005,shish2009,kop2010,geng2016}).
At the present time, there are no effective numerical methods that allow directly
(without a preliminary nonlinear point transformation of the given equation)
to consider such hypersingular problems even for moderately small $\varepsilon\simeq\frac 1{50}$
and, especially, for smaller values of $\varepsilon$.

\begin{remark}
A qualitatively similar situation occurs for a hypersingular boundary-value problem described by
the equation $\varepsilon y^{\prime\prime}_{xx}+(y^\prime_x)^2-1=0$
with the boundary conditions \eqref{e01a}.
\end{remark}

\begin{remark}
Equation~\eqref{e01} describes a family of the traveling wave exact solutions
of the potential Burgers equation $u_t+u_z^2=\varepsilon u_{zz}$ \cite{Pol_Zai2012},
in which we have to set $u=-py(x)$, where $x=z-qt$.
\end{remark}

\begin{remark}
Differential equation \eqref{e01} do not satisfy the Painlev\'e test \cite{cont1999,kudr2010}.
\end{remark}

%%%%%%%%%%%%%%%%%%%%
\subsection{Two classes of hypersingular problems admitting linearization}\label{ss:2.2}

\textit{Problem 2.} We consider the nonlinear boundary-value problem with a small parameter
\begin{align}
\varepsilon y^{\prime\prime}_{xx}+p(y^\prime_x)^2+q(x)y^\prime_x+r(x)=0;\quad \ y(0)=\alpha,\quad \ y(1)=\beta,
\label{204}
\end{align}
which generalizes the problem \eqref{e01}--\eqref{e01a} and contains two arbitrary functions
$q(x)$, $r(x)$ and the constant $p$.

The substitution $y=(\varepsilon/p)\ln u$
reduces the problem \eqref{204} to the linear problem:
\begin{align}
\varepsilon^2u^{\prime\prime}_{xx}+\varepsilon q(x)u^\prime_x+pr(x)u=0;\quad \ u(0)=A,\quad \ u(1)=B,
\label{206}
\end{align}
where $A=e^{\alpha p/\varepsilon}$ and $B=e^{\beta p/\varepsilon}$.

If $q(x)=\text{const}$ and  $r(x)=\text{const}$, then we can obtain the exact
solution of the problem~\eqref{206}.
If $q$ and/or $r$ are not constants, then to obtain an approximate solution as $\varepsilon\to 0$
we can use the method of matched asymptotic expansions
\cite{eck1973,lag1988,mal1991,ker1996,nay2000,Pol_Zai2003}.

\begin{remark}
The problem \eqref{204}, as well as the problem \eqref{e01}--\eqref{e01a},
does not allow a direct application of the method of matched asymptotic expansions
(without a preliminary nonlinear point transformation of the equation under consideration),
since the second term of this equation remains dominant over the first term
with any change in the scale of the independent variable by means of the substitution
$x=k\bar x$, where $k=k(\varepsilon)$.
\end{remark}

\textit{Problem 3.}
The nonlinear boundary-value problem with a small parameter
\begin{align}
\varepsilon y^{\prime\prime}_{xx}+p(y)(y^\prime_x)^2+q(x)y^\prime_x=0;\quad \ y(0)=\alpha,\quad \ y(1)=\beta,
\label{211aa}
\end{align}
containing two arbitrary functions $p(y)$ and $q(x)$, with the substitution
\begin{align}
u=P(y),\quad \ P(y)=\int^y_0\exp\BL[\frac 1{\varepsilon}\int^z_0p(\xi)\,d\xi\BR]dz,
\label{211ab}
\end{align}
reduces to the linear problem
\begin{align}
\varepsilon u^{\prime\prime}_{xx}+q(x)u^\prime_x=0;\quad \ u(0)=P(\alpha),\quad u(1)=P(\beta),
\label{211ac}
\end{align}
that is easily integrable.

In particular, for $q=\text{const}$, by using the transformation~\eqref{211ab} and equation~\eqref{211ac},
we obtain the solution of the problem~\eqref{211aa} in an implicit form
\begin{align}
P(y)=\frac{P(\beta)-P(\alpha)e^{-q/\varepsilon}}{1-e^{-q/\varepsilon}}+\frac{P(\alpha)-P(\beta)}{1-e^{-q/\varepsilon}}e^{-qx/\varepsilon}.
\label{211ad}
\end{align}
%ïðîâåðåíî!
Differentiating~\eqref{211ad}, we find the derivative on the left boundary:
\begin{align}
y^\prime_x\br|_{x=0}=
\frac q{\varepsilon}\frac{P(\beta)-P(\alpha)}{P'_y(\alpha)(1-e^{-q/\varepsilon})}
=\frac q{\varepsilon}\frac{P(\beta)-P(\alpha)}{(1-e^{-q/\varepsilon})}\exp\BL[-\frac 1{\varepsilon}\int^\alpha_0p(\xi)\,d\xi\BR].
\label{211ae}
\end{align}
%In the particular case $\alpha=0$, the formula~\eqref{211ae} is simplified and takes the form
%\begin{align}
%y^\prime_x\br|_{x=0}=
%\frac q{\varepsilon}\frac{P(\beta)}{1-e^{-q/\varepsilon}}.
%\label{211af}
%\end{align}

\begin{remark}
In the problem~\eqref{211aa} and the formulas~\eqref{211ab},~\eqref{211ad}--\eqref{211ae},
the functions $p(y)$ and $q(x)$ may also depend on the small parameter $\varepsilon$.
\end{remark}
\medskip

\textit{Example.}
We consider the problem~\eqref{211aa} for $p(y)=e^y+\varepsilon$, $q(x)=1$, $\beta>0$, $\alpha=0$ ($\varepsilon\to 0$).
Substituting these functions and constants in~\eqref{211ae}, we obtain a very large derivative on the boundary
\begin{align}
y^\prime_x\br|_{x=0}=
\frac{1}{1-e^{-1/\varepsilon}}\BL[\exp\BL(\frac{e^\beta-1}\varepsilon\BR)-1\BR]\simeq \exp\BL(\frac{e^\beta-1}\varepsilon\BR).
\label{211ag}
\end{align}
In this problem, the hypersingular boundary layer is more thin than in the problem \eqref{e01}--\eqref{e01a} for $p=q=1$.

%%%%%%%%%%%%%%%%%%%%%%%%% 3

\section{Hypersingular boundary-value problems for partial differential equations}\label{s:3}

\subsection{Example of a hypersingular boundary-value problem. Exact solution}\label{ss:3.1}

\textit{Problem 4.}
We consider the initial-boundary value problem with a small parameter for the
potential Burgers equation \cite{Pol_Zai2012}:
\begin{align}
w_t=\varepsilon w_{xx}+w_x^2\quad \ (t>0, \ x>0);\quad \ w(0,t)=\alpha,\quad w(x,0)=\beta.
\label{300}
\end{align}
This problem has a self-similar solution of the form
\begin{align}
w=W(z),\quad \ z=xt^{-1/2},
\label{301}
\end{align}
where the function $W(z)$ is described by the ordinary differential equation
\begin{align}
\varepsilon W''_{zz}+(W'_z)^2+\tfrac12zW'_z=0;\quad \ W(0)=\alpha,\quad \ W(\infty)=\beta,
\label{302}
\end{align}
which belongs to the class of equations~\eqref{204} and admits linearization
by the substitution $W=\varepsilon \ln u$.

The exact solution of the problem~\eqref{302} is expressed in terms of the error function (probability integral):
\begin{align}
{}\!\!\!W=\varepsilon\ln\BL[(e^{\beta/\varepsilon}-e^{\alpha/\varepsilon})\text{erf}\BL(\frac z{2\sqrt{\varepsilon}}\BR)
+e^{\alpha/\varepsilon}\BR], \ \ \
\text{erf}(\xi)=\frac 2{\sqrt{\pi}}\int^\xi_0\exp\bl(-\zeta^2\br)\,d\zeta.
\label{303}
\end{align}
Calculating the derivative on the left boundary, we have
\begin{align}
W'_z|_{z=0}=\sqrt{\varepsilon/\pi}\bl[e^{(\beta-\alpha)/\varepsilon}-1\br].
\label{304}
\end{align}
For $\beta>\alpha$ and $\varepsilon\to 0$, the boundary layer near $x=0$
has a hypersingular character, and for $\beta\ge\alpha$ there is no boundary layer.

%%%%%%%%%%%%%%%%%%%
\subsection{Two classes of hypersingular problems admitting linearization}\label{ss:3.2}

\textit{Problem 5.}
We consider a stationary hypersingular boundary-value problem with a small parameter
for the partial differential equation of elliptic type:
\begin{align}
\varepsilon \Delta w+p|\nabla w|^2+\textbf{q}(\textbf{x})\cdot\nabla w+r(\textbf{x})=0,
\quad \text{$\textbf{x}\in V$};\quad w|_{\partial V_1}=\alpha,\quad w|_{\partial V_2}=\beta.
\label{400}
\end{align}
Here $\Delta$ is the Laplace operator; $\nabla$ is the Hamilton operator (gradient operator);
$\alpha$, $\beta$, $p$ are constants, $\textbf{x}=(x_1,\dots,x_n)$; $\textbf{q}=\textbf{q}(\textbf{x})$
and $r=r(\textbf{x})\not\equiv 0$ are given functions, $V$ is an open domain bounded by the surfaces
$\partial V_1$ and $\partial V_2$.

The substitution $w=(\varepsilon/p)\ln u$ reduces the problem~\eqref{400} to the linear problem:
\begin{align}
\varepsilon^2 \Delta u+\varepsilon\textbf{q}(\textbf{x})\cdot\nabla u+pr(\textbf{x})u=0,
\quad\text{$\textbf{x}\in V$};\quad \ w|_{\partial V_1}=A,\quad w|_{\partial V_2}=B,
\label{401}
\end{align}
where $A=e^{\alpha p/\varepsilon}$ and $B=e^{\beta p/\varepsilon}$.

\begin{remark}
The partial differential equation of parabolic type
\begin{equation*}
w_t=\varepsilon \Delta w+p|\nabla w|^2+\textbf{q}(\textbf{x},t)\cdot\nabla w+r(\textbf{x},t)
\end{equation*}
is also linearized by the substitution $w=(\varepsilon/p)\ln u$.
\end{remark}

\textit{Problem 6.}
The nonlinear boundary-value problem
\begin{align*}
\varepsilon \Delta w+p(w)|\nabla w|^2+\textbf{q}(\textbf{x},w)\cdot\nabla w=0,\quad\text{$\textbf{x}\in V$};\quad \
w\br|_{\partial V_1}=\alpha,\quad w\br|_{\partial V_2}=\beta,
%\label{501aa}
\end{align*}
with the substitution $u=P(w)$, where the function $P(w)$ is defined in \eqref{211ab},
%\begin{align*}
%u=P(w),\quad \ P(w)=\int^w_0\exp\BL[\frac 1{\varepsilon}\int^z_0p(\xi)\,d\xi\BR]dz,
%\label{501ab}
%\end{align*}
reduces to a linear problem.

%%%%%%%%%%%%%%%%% 4

\section{Brief conclusions}\label{s:4}

For the first time, hypersingular nonlinear boundary-value problems
with a small parameter $\varepsilon$ are described,
in which hyperfine boundary layers arise and a solution on the boundary of the layer
has a very large derivative of the order of $e^{1/\varepsilon}$ and more
(in standard problems with boundary layers, the derivative on the boundary is of the order of $\varepsilon^{-1}$ or less).
Exact solutions of some hypersingular boundary-value problems for ODEs
and PDEs are obtained.

%The authors thank A. V. Aksenov for his interest in the work and for a useful discussion.

\end{document}